\documentclass[shortnote,twocolumn]{jpsj3}
\usepackage{txfonts}
\usepackage{bm}

\title{Tit-for-Tat Strategy as a Deformed Zero-Determinant Strategy in Repeated Games}

\author{Masahiko Ueda$^1$\thanks{m.ueda@yamaguchi-u.ac.jp}}
\inst{$^1$Graduate School of Sciences and Technology for Innovation, Yamaguchi University, Yamaguchi 753-8511, Japan} %\\

\abst{We introduce the concept of deformed zero-determinant strategies in repeated games.
We then show that the Tit-for-Tat strategy in the repeated prisoner's dilemma game is a deformed zero-determinant strategy, which unilaterally equalizes the probability distribution functions of payoffs of two players.}

\begin{document}
\maketitle

The Tit-for-Tat (TFT) strategy is a strategy in the iterated prisoner's dilemma game which chooses the opponent's previous action \cite{RCO1965}.
It has been known that TFT forms the Nash equilibrium, where both players are cooperative.
Although the strategy is simple, it attained the highest average score in computer tournaments \cite{AxeHam1981}.
Its role in evolutionary game theory has substantially been investigated \cite{NowSig1992,NowSig1993,IFN2007}.
Furthermore, it was recently pointed out that TFT is contained in a class of memory-one strategies called zero-determinant (ZD) strategies \cite{PreDys2012}, which unilaterally enforce linear relations between average payoffs of players.
Although TFT is not robust against implementation errors, TFT was used to construct a longer-memory strategy which is successful even if implementation errors exist \cite{YBC2017}.
It has also been known that TFT cannot be beaten in several situations \cite{DOS2014}.

In this paper, we introduce the concept of deformed zero-determinant (DZD) strategies in repeated games.
We then show that the TFT strategy is a DZD strategy which unilaterally enforces linear relations between all moments of payoffs of two players, which implies that the probability distribution functions of payoffs of two players are equal to each other.
This result provides a fresh perspective on the TFT strategy.

We consider the iterated prisoner's dilemma game \cite{PreDys2012}.
There are two players ($1$ and $2$) in the game.
Each player takes cooperation (described as $C$) or defection (described as $D$) in a one-shot game.
The action of player $a$ is written as $\sigma_a \in \{ C, D \}$.
We collectively write $\bm{\sigma}:=\left( \sigma_1, \sigma_2 \right)$.
The payoff of player $a\in \{ 1, 2 \}$ when the state is $\bm{\sigma}$ is described as $s_a\left( \bm{\sigma} \right)$.
The payoffs in the prisoner's dilemma game are defined as
$\bm{s}_1:=\left( s_1 \left( C, C \right), s_1 \left( C, D \right), s_1 \left( D, C \right), s_1 \left( D, D \right) \right)=(R,S,T,P)$ and $\bm{s}_2:=\left( s_2 \left( C, C \right), s_2 \left( C, D \right), s_2 \left( D, C \right), s_2 \left( D, D \right) \right)=(R,T,S,P)$
with $T>R>P>S$ and $2R>T+S$.
The memory-one strategy of player $a$ is described as the conditional probability $T_a\left( \sigma_a | \bm{\sigma}^{\prime} \right)$ of taking action $\sigma_a$ when the state in the previous round is $\bm{\sigma}^{\prime}$.
Then, the time evolution of this system is described as the Markov chain
\begin{eqnarray}
 P\left( \bm{\sigma}, t+1 \right) &=& \sum_{\bm{\sigma}^{\prime}} T\left( \bm{\sigma} | \bm{\sigma}^{\prime} \right) P\left( \bm{\sigma}^{\prime}, t \right)
\end{eqnarray}
with the transition probability
\begin{eqnarray}
 T\left( \bm{\sigma} | \bm{\sigma}^{\prime} \right) &:=& \prod_{a=1}^2 T_a\left( \sigma_a | \bm{\sigma}^{\prime} \right),
\end{eqnarray}
where $P\left( \bm{\sigma}^{\prime}, t \right)$ is the probability distribution of a state $\bm{\sigma}^{\prime}$ at time $t$.
We consider the case that there is no discounting of future payoffs.

We now introduce the concept of deformed zero-determinant (DZD) strategies.
The original ZD strategies of player $a$ are strategies which can be written in the form
\begin{eqnarray}
 \sum_{\sigma_a} c_{\sigma_a} \hat{T}_a\left( \sigma_a | \bm{\sigma}^{\prime} \right) &=& \sum_{b=0}^2 \alpha_{b} s_{b} \left( \bm{\sigma}^{\prime} \right) 
\end{eqnarray}
with some coefficients $\left\{ \alpha_{b} \right\}$ and $\left\{ c_{\sigma_a} \right\}$, where we have defined
\begin{eqnarray}
 \hat{T}_a\left( \sigma_a | \bm{\sigma}^{\prime} \right) &:=& T_a\left( \sigma_a | \bm{\sigma}^{\prime} \right) -  \delta_{\sigma_a, \sigma^{\prime}_a}
\end{eqnarray}
and $s_0\left( \bm{\sigma} \right):=1$.
The term $\delta_{\sigma, \sigma^\prime}$ is the Kronecker delta.
(It should be noted that $\sum_{\sigma_a} \hat{T}_a\left( \sigma_a | \bm{\sigma}^{\prime} \right)=0$ for any $\bm{\sigma}^{\prime}$, due to the normalization condition of probability.)
The quantity $\hat{T}_a\left( \sigma_a | \bm{\sigma}^{\prime} \right)$ is called as a Press-Dyson vector.
Since the average of $\hat{T}_a\left( \sigma_a | \bm{\sigma}^{\prime} \right)$ with respect to the stationary distribution $P^{(\mathrm{st})}\left( \bm{\sigma}^{\prime} \right)$ is zero \cite{Aki2012,UedTan2020}
\begin{eqnarray}
 0 &=& \sum_{\bm{\sigma}^{\prime}} \hat{T}_a\left( \sigma_a | \bm{\sigma}^{\prime} \right) P^{(\mathrm{st})}\left( \bm{\sigma}^{\prime} \right),
\end{eqnarray}
ZD strategies unilaterally enforce a linear relation between average payoffs:
\begin{eqnarray}
 0 &=& \sum_{b=0}^2 \alpha_{b} \left\langle s_{b}  \right\rangle^{(\mathrm{st})},
\end{eqnarray}
where $\left\langle \cdots \right\rangle^{(\mathrm{st})}$ represents the average with respect to the stationary distribution $P^{(\mathrm{st})}\left( \bm{\sigma}^{\prime} \right)$.
Now, we introduce DZD strategies as ones satisfying
\begin{eqnarray}
 \sum_{\sigma_a} c_{\sigma_a} \hat{T}_a\left( \sigma_a | \bm{\sigma}^{\prime} \right) &=& \sum_{k_1=0}^\infty \sum_{k_2=0}^\infty \alpha^{(k_1, k_2)} s_{1} \left( \bm{\sigma}^{\prime} \right)^{k_1} s_{2} \left( \bm{\sigma}^{\prime} \right)^{k_2}
\end{eqnarray}
with some coefficients $\left\{ \alpha^{(k_1, k_2)} \right\}$ and $\left\{ c_{\sigma_a} \right\}$.
Then, DZD strategies unilaterally enforce a linear relation between moments of payoffs:
\begin{eqnarray}
 0 &=& \sum_{k_1=0}^\infty \sum_{k_2=0}^\infty \alpha^{(k_1, k_2)} \left\langle s_{1}^{k_1} s_{2}^{k_2} \right\rangle^{(\mathrm{st})}
\end{eqnarray}
Even if payoff vectors $\bm{s}_1$ and $\bm{s}_2$ and the vector of all ones $\bm{1}:=(1, 1, 1, 1)$ do not form a basis that spans the space of all Press-Dyson vectors, this extension of the basis generally enables any Press-Dyson vectors to be represented by the basis vectors.
Although we introduced the concept of DZD strategies for the repeated prisoner's dilemma game, extension to general multi-player multi-action games is straightforward.

Concretely, we consider the TFT strategy of player $1$:
\begin{eqnarray}
 \bm{T}_1 (C) &:=& \left(
\begin{array}{c}
T_1 \left( C | C, C \right)  \\
T_1 \left( C | C, D \right)   \\
T_1 \left( C | D, C \right)   \\
T_1 \left( C | D, D \right) 
\end{array}
\right)
=
\left(
\begin{array}{c}
1  \\
0  \\
1   \\
0
\end{array}
\right).
\end{eqnarray}
(Although the TFT strategy is not a ZD strategy in general under observation errors \cite{MamIch2019}, it is a ZD strategy enforcing $0 = \left\langle s_{1} \right\rangle^{(\mathrm{st})} - \left\langle s_{2} \right\rangle^{(\mathrm{st})}$ when there are no errors \cite{PreDys2012}.)
Then, her Press-Dyson vector is written as
\begin{eqnarray}
 \hat{\bm{T}}_1 (C) &:=& \left(
\begin{array}{c}
\hat{T}_1 \left( C | C, C \right)  \\
\hat{T}_1 \left( C | C, D \right)   \\
\hat{T}_1 \left( C | D, C \right)   \\
\hat{T}_1 \left( C | D, D \right) 
\end{array}
\right)
=
\left(
\begin{array}{c}
0  \\
-1  \\
1   \\
0
\end{array}
\right)
\end{eqnarray}
and we find
\begin{eqnarray}
  \hat{\bm{T}}_1 (C) &=& \frac{1}{T^k-S^k} \left[ \bm{s}_1^k - \bm{s}_2^k \right]
\end{eqnarray}
for arbitrary $k\geq 1$, where we have introduced the notation $\bm{s}_a^k:=\left( s_a \left( C, C \right)^k, s_a \left( C, D \right)^k, s_a \left( D, C \right)^k, s_a \left( D, D \right)^k \right)$.
Therefore, the TFT strategy is contained in DZD strategies, and we obtain a linear relation
\begin{eqnarray}
  0 &=& \left\langle s_{1}^{k} \right\rangle^{(\mathrm{st})} - \left\langle s_{2}^{k} \right\rangle^{(\mathrm{st})}.
  \label{eq:linear_DZDS_TFT}
\end{eqnarray}
In other words, the TFT strategy unilaterally enforces linear relations between all moments of payoffs of two players.
Although the case of $k=1$ was known in Ref. \cite{PreDys2012}, we find that Eq. (\ref{eq:linear_DZDS_TFT}) holds for any $k\geq 1$.
We remark that the strategy of player $2$ is arbitrary.

From another point of view, when we introduce the quantity $\bm{\Phi}_a(h):=\left( e^{h s_a \left( C, C \right)}, e^{h s_a \left( C, D \right)}, e^{h s_a \left( D, C \right)}, e^{h s_a \left( D, D \right)} \right)$, the TFT strategy satisfies
\begin{eqnarray}
  \hat{\bm{T}}_1 (C) &=& \frac{1}{e^{hT}-e^{hS}} \left[ \bm{\Phi}_1(h) - \bm{\Phi}_2(h) \right]
\end{eqnarray}
for $h\neq 0$.
Therefore, we obtain a linear relation
\begin{eqnarray}
  0 &=& \left\langle e^{h s_{1}} \right\rangle^{(\mathrm{st})} - \left\langle e^{h s_{2}} \right\rangle^{(\mathrm{st})},
\end{eqnarray}
which means that the moment generating functions of payoffs of two players are equal to each other under the TFT strategy.
This is the main result of this paper.
Because the equality of two moment generating functions implies the equality of two probability distribution functions, we conclude that TFT unilaterally enforces equality of the probability distributions of payoffs of two players.

We finally remark that, although all memory-one strategies are not necessarily ZD strategies, they are DZD strategies in general.
For instance, it is known that the win-stay lose-shift (WSLS) strategy \cite{NowSig1993}
\begin{eqnarray}
 \bm{T}_1 (C) &=&
\left(
\begin{array}{c}
1  \\
0  \\
0   \\
1
\end{array}
\right)
\end{eqnarray}
is not contained in the class of ZD strategies.
However, it is generally contained in the class of DZD strategies.
For example, its Press-Dyson vector is described as
\begin{eqnarray}
  \hat{\bm{T}}_1 (C) &=& \alpha_1 \bm{s}_1 + \alpha_2 \bm{s}_2 + \alpha_{1,2} \bm{r}_{1,2} + \gamma \bm{1}
\end{eqnarray}
with the appropriate coefficients $\left( \alpha_1,  \alpha_2,  \alpha_{1,2}, \gamma \right)$, where we have introduced the notation
\begin{eqnarray}
 \bm{r}_{1,2} &:=&
\left(
\begin{array}{c}
s_1 \left( C, C \right)s_2 \left( C, C \right)  \\
s_1 \left( C, D \right)s_2 \left( C, D \right)  \\
s_1 \left( D, C \right)s_2 \left( D, C \right)   \\
s_1 \left( D, D \right)s_2 \left( D, D \right)
\end{array}
\right),
\end{eqnarray}
because the dimension of the space of Press-Dyson vectors $\hat{\bm{T}}_1 (C)$ is four and the vectors $\left( \bm{s}_1, \bm{s}_2, \bm{r}_{1,2}, \bm{1} \right)$ are linearly independent in general.
Then, the WSLS strategy unilaterally enforces
\begin{eqnarray}
  0 &=& \alpha_1 \left\langle s_1 \right\rangle^{(\mathrm{st})} + \alpha_2 \left\langle s_2 \right\rangle^{(\mathrm{st})} + \alpha_{1,2} \left\langle s_1 s_2 \right\rangle^{(\mathrm{st})} + \gamma.
\end{eqnarray}
However, this linear relation is nonsense because the coefficients $\left( \alpha_1,  \alpha_2,  \alpha_{1,2}, \gamma \right)$ depend on the concrete values of payoffs $(R,S,T,P)$.
Moreover, similarly as the case of the TFT strategy, DZD strategies generally enforce multiple linear relations between moments of payoffs simultaneously.

In this paper, we introduced the concept of DZD strategies in repeated games.
We then proved that the TFT strategy is a DZD strategy which unilaterally equalizes the moment generating functions of payoffs of two players.
We believe that this result deepens the understanding of the TFT strategy.
Constructing useful examples of DZD strategies is an important future problem.

\begin{acknowledgment}

%\acknowledgment

This study was supported by JSPS KAKENHI Grant Number JP20K19884.

\end{acknowledgment}

\bibliographystyle{jpsj}
\bibliography{deformedZDS}

\end{document}